\documentclass[a4paper,12pt]{article}
\usepackage[dvips]{graphicx}
\usepackage{latexsym}
\usepackage{amsfonts}
\usepackage{amssymb}
\usepackage{amsmath}
\usepackage{amsthm}
\usepackage{pstricks}

\newtheorem{theorem}{Theorem}[section]
\newtheorem{lemma}[theorem]{Lemma}
\newtheorem{definition}[theorem]{Definition}

\def\E{\mathbb{E}}
\def\F{ {\mathcal{F} }}

\def\J{ {\mathcal{J} }}

\def\N{ {\mathbb{N} }}
\def\P{\mathbb{P}}

\def\Z{ {\mathbb{Z} }}

\begin{document}

\title{Universal Relationships in Measures of Unpredictability}

\author{Finn Macleod\thanks{MACSI, University of Limerick, Ireland}, Alexei Pokrovskii$^{\dagger}$, Dmitrii Rachinskii\thanks{Department of Applied
Mathematics, University College Cork, Ireland; Institute for Information Transmission Problems,
Russian Academy of Sciences, {\em on leave}}}

\date{}

\maketitle
\abstract

The predictability of a sequence is defined as the asymptotic performance of the best performing predictor in a given class. The value of the predictability of a sequence will in general depend on the choice of this predictor class. The existence of universal properties of predictability is demonstrated by looking at relationships between different sequences - these relationships hold for any class of predictors satisfying a certain set of axioms.


\medskip
{\bf Mathematical Subject Classification:} {\bf ??}

\section{Introduction}
\label{s:randomness}

How \emph{predictable} is a given sequence of digits? Certainly some sequences,
\begin{equation*}
0000000000\ldots
\end{equation*}
seem more predictable than others,
\begin{equation*}
0110101011\ldots ,
\end{equation*}
in the same way as some sequences appear more random than others. However, characterising predictability is a question that is distinct from notions of randomness arising in the more well known areas of probability theory and Kolmogorov complexity. One can consider three different meanings of the word random:

\begin{enumerate}
\item In probability theory, a random sequence is as a result of a `random selection' from some set - the randomness is a property of the measure on the set.

\item Descriptive, or Kolmogorov complexity. The Kolmogorov complexity of a sequence is the length of the shortest method for describing that sequence. A sequence which has no method of description shorter than itself is considered random.

\item Predictability. A sequence is random if it is difficult to predict.
\end{enumerate}

The links between randomness in probability theory and that of Kolmogorov complexity are well known. They arise via Shannon entropy, for example, with high probability, sequences chosen from a set will have Kolmogorov complexity close to the Shannon entropy. See \cite{CT91} or \cite{livitanyi} for a brief introduction to these ideas.

Bounds are also known which link the Kolmogorov complexity to our notion of predictability (defined below). However the two quantities are distinct, and there exist sequences with the same Kolmogorov complexity and different predictability, and vice versa \cite{D73, FMG92}.

The definition of predictability we discuss was first introduced in \cite{P89}. It arose independently in \cite{FMG92} using a specific predictor class. We use the binary setting: $\{0,1\}^{\infty}$ denotes the space of all binary sequences $a=a_0a_1a_2a_3\ldots$

\begin{definition}
A binary predictor is any mapping between two infinite binary sequences
\begin{equation*}
f:\{0,1\}^{\infty} \to \{0,1\}^{\infty}
\end{equation*}
with the property of \emph{causality}; that is, $(f(a))_0$ is the same for all $a\in \{0,1\}^{\infty}$ and for each $n\ge 1$ given $a = a_0a_1a_2\ldots$ and  $b = b_0b_1b_2\ldots \in \{0,1\}^{\infty}$ with $a_i = b_i$ for $i = 0, \ldots n-1$, then
\begin{equation*}
(f(a))_n = (f(b))_n.
\end{equation*}
\end{definition}
We equip a class of predictors with a \textit{hierarchy}.

\begin{definition}
A predictor hierarchy on $\F$ is a set of increasing sets of predictors, $\F_1,\F_2\ldots$, with $\F_i \subset \F_{i+1}$ and $\bigcup_{i=1}^{\infty}\F_i = \F$.
\end{definition}

We now define predictability as the accuracy of the best performing predictor in a given class. These classes can be infinite - for example that of finite state automata, or all computable prediction strategies (see Section \ref{d:classes}). Thus we approach any value of predictability asymptotically, and use the idea of a hierarchy to enable this. In the latter case, we note that predictability, like Kolmogorov complexity, will not be a computable quantity.
\begin{definition}
The predictability $I(a; \F)$ of a sequence $a$ with respect to a predictor hierarchy $\F$ is
\begin{equation}\label{IaF}
I(a;\F) = \lim_{m \to \infty} \limsup_{n \to \infty} \min_{f \in
\F_m} \frac{1}{n}\sum_{i=0}^{n-1}((f(a))_i\oplus a_i)
\end{equation}
where $\oplus$ denotes summation mod 2.
\end{definition}

One can show that the predictability is independent of the hierarchy chosen, but it is dependant on the class of predictors. As an example, consider the binary expansion of $\pi$. It can be predicted perfectly by an algorithm which generates the digits of $\pi$, but no finite state machine has the unbounded memory to do this, and thus will accrue errors. Thus the predictability of $\pi$ with respect to the two hierarchies of finite state automata and computable prediction strategies will differ. That predictability is independent of the hierarchy chosen follows from the definitions. We attach details in Appendix 3.

However, we might still believe that some \emph{operations} on sequences universally increase or decrease predictability, irrespective of predictor class. Consider a sequence $a = a_0a_1a_2a_3 \ldots$, and form the new sequence
\begin{equation*}
\begin{tabular}{cccccccccc}
$b=S(a) =$ & $a_0 \oplus a_1$ & $a_3\oplus a_4$ &  $a_6 \oplus a_7$
& $a_9\oplus a_{10}$& $\ldots$
\end{tabular}
\end{equation*}
The digits are mixed together, and given $b$, we can not determine the sequence $a$. In general, one would expect this kind of operation to make a sequence less predictable. But it is also possible that the sequence $a$ is more predictable. For example, take $a$ with $a_{3i}=a_{3i+1}=1$  and allow only $a_{3i+2}$ to vary.  Under the operation $S$, we will obtain a perfectly predictable, constant sequence.

We claim that if one has a sequence which becomes more predictable under the operation $S$, then that says something about the structure of $a$; the structure of $a$ is somehow linked to the structure of the operation $S$. This is the idea behind our central result. Either:
\begin{enumerate}
\item Certain simple operations on a sequence will cause a sequence to be more difficult to predict, or
\item There exists a subsequence of $a$ which is easier to predict than $a$.
\end{enumerate}

We establish this theorem with the use of some general axioms about a predictor hierarchy.

We will say that a sequence is independent if there is no rule (in terms of predictors from the class ${\mathcal F}$)
for selecting a subsequence with a different value of predictability. Thus for independent sequences the above theorem simplifies. We will prove a corollary which enables comparisons with analagous ideas in probability theory.

\section{Existence of all values of predictability}

We assume that the class $\F$ contains the constant mappings $\phi^0, \phi^1$ defined by $(\phi^0(a))_n = 0,  (\phi^1(a))_n = 1$. Therefore for any $a \in \{0,1\}^{\infty}$, $ I(a)\in [0,1/2]$. Then we can show the following.

\begin{theorem}\label{tt2}
For any $I_0 \in [0,1/2]$ there exist sequences $a \in \{0,1\}^{\infty}$ which satisfy $I(a;\F) = I_0$.
\end{theorem}

The proof of this theorem is relegated to Appendix 1. We conjecture that though there exist sequences taking all values of unpredictability between $0$ and $1/2$, almost all (in the probabilistic sense) will have unpredictability $1/2$. Indeed we can imagine large deviations type arguments where, if we consider any restricted set consisting of sequences taking unpredictability values in $[a,b]$ with $a<b$, then almost all the sequences in that set will take the larger unpredictability value $b$.

{\bf D: Alexei, care to comment on the above paragraph. Can we state this fact, not conjecture?}

We now introduce the axioms we require to establish our central result.

\section{Axioms of Predictor hierarchies}

These axioms are the weakest set of assumptions required to prove our theorem. We will sometimes write $fa$ rather than $f(a)$, when it is clear that the predictor $f$ is acting on $a$.

We first define the following operations on sequences.

\begin{definition}
We define two operations:
\begin{enumerate}
\item The extraction of subsequences. For $\nu = 0,1,2$, define $P^{\nu}: \{0,1\}^\infty \to \{0,1\}^\infty$ with
$(P^{0}a)_i = a_{3i }$, $(P^{1}a)_i = a_{3i +2}$, $(P^{2}a)_i =
a_{3i +1}$.
\item Summation of subsequences. For $\nu = 1,2$ define $S^{\nu}: \{0,1\}^\infty \to \{0,1\}^\infty$ with
$(S^{1}a)_i= a_{3i}\oplus a_{3i + 2}$, $(S^{2}a)_i= a_{3i+1}\oplus a_{3i + 2}$.
\end{enumerate}
\end{definition}

For example, for any sequence $a=a_0a_1a_2a_3...$
$$
P^1(a)=a_2\ a_5\ ...
$$
$$
S^1(a)=(a_0\oplus a_2)\ \ (a_3\oplus a_5)\ ...
$$

We now introduce a method for selecting subsequences from a sequence using a predictor.

\begin{definition}
The subsequence selected from $a$ by predictor $f$, $f_* a$, is a sequence $b \in \{0,1\}^{\infty}$,
defined by $b_l = a_{i(l)}$, where $i(l)$ specifies the $l$th-index for which $(fa)_i = 1$ holds.
\end{definition}

Whenever $f$ takes the value 1, that digit is added to the
subsequence. For example, if $f$ is periodic predictor, predicting
$0011$ periodically, independent of input, then if
$a=a_0a_1a_2a_3\ldots$
\begin{equation*}
f_*a = a_2a_3a_6a_7a_{10}a_{11}\ldots
\end{equation*}

We now state the axioms we require.

\begin{description}
\item[Axiom 1 (Summation).] For any $f^0, f^1 \in \F$, $\F$ also contains the mapping $f = f^0\oplus f^1$ defined by
\begin{equation*}
(fa)_i = (f^0a)_i \oplus (f^1a)_i.
\end{equation*}
\item[Axiom 2 (Interleaving).] For any $f^0, f^1, f^2 \in \F$, $\F$ also contains the mapping $f$ defined by the relation
\begin{equation*}
(fa)_{3i-\nu} = (f^{\nu}a)_{3i- \nu}
\end{equation*}
for $\nu = 0,1,2$. Equivalently,
$$
P^0fa=P^0f^0a,\quad P^1fa=P^1f^1a,\quad P^2fa=P^2f^2a.
$$
\item[Axiom 3 (Subsequences).] For any $f \in \F$ , the class $\F$ also contains at least one mapping, $f^1$, which satisfies
\begin{eqnarray*}
P^1f^1 a = fS^1a,
\end{eqnarray*}
at least one mapping, $f^2$, which satisfies
\begin{eqnarray*}
P^1f^2 a = fS^2a,
\end{eqnarray*}
at least one mapping, $g^1$, which satisfies
\begin{eqnarray*}
P^2g^1 a = fS^1a,
\end{eqnarray*}
and at least one mapping, $g^2$, which satisfies
\begin{eqnarray*}
P^2 g^2a = fS^2a.
\end{eqnarray*}
Similarly, for any $f \in \F$, $\F$ also contains at least one mapping, $h^0$, which satisfies:
\begin{eqnarray*}
P^0h^0a = fP^0a,
\end{eqnarray*}
at least one mapping, $h^1$, which satisfies
\begin{eqnarray*}
P^1h^1a = fP^1a
\end{eqnarray*}
and at least one mapping, $h^2$, which satisfies
\begin{eqnarray*}
P^2h^2a = fP^2a.
\end{eqnarray*}
\item[Axiom 4 (Switching).] For any $f^0, f^1, f^2 \in \F$, $\F$ also contains the mapping $f$ specified by
\begin{equation*}
(fa)_i = \left\{ \begin{array}{ll}
         (f^1a)_i & \mbox{if $(f^0a)_i = 0$},\\
        (f^2b)_{l(i)} & \mbox{if $(f^0a)_i = 1$},\end{array} \right.
\end{equation*}
where sequence $b$ is defined by $b = f^0_*a$; $l(i)$ is the number
of indices $j$ which satisfy the relations $j<i, (f^0a)_j = 1$.
At each point where $(f^0a)_i = 1$, this indexing system selects
sequentially elements from the sequence $(f^2b)_0, (f^2b)_1,
 \ldots$, which is what we require. \label{axiom4}
\end{description}

We will assume Axioms 1--4 to hold. We will also assume that the class $\F$ contains the constant predictors
$\phi^0,\phi^1$ and the simple predictors
\begin{equation}\label{psipsi}
(\psi^1 a)_j=a_{j-2},\qquad (\psi^2a)_j=a_{j-1}.
\end{equation}

\section{Examples of predictor hierarchies}
We have two examples in mind when considering classes of predictors which satisfy the above axioms:
\begin{enumerate}
\item Finite state automata.

\item The class of all computable predictors based on Turing machines.
\end{enumerate}
We prove Axioms 1-4 for the class of all finite state automata and sketch the proof
for the class of computable predictors in Appendix 2.

Notably, the class of Markov predictors does not satisfy Axiom 4. Axiom 4 requires that the predictors have the capacity to base their predictions upon events arbitrarily far back in the past. Markov predictors do not have this property - they make their predictions based purely on a finite window of time. Other potential candidates for predictor classes satisfying our axioms can be surmised from language theory: for example, pushdown automata or linear bounded automata (these both contain finite state automata as a subset).

\label{d:classes}

\section{Unpredictability relationships of sequences}

\begin{definition}
The fraction of the first $n$ terms of a sequence $a$ which take the value 1 is given by
\begin{equation*}
E(a;n)= \frac{1}{n}\sum_{i=0}^{n-1} a_i.
\end{equation*}
\end{definition}
We are now in a position to prove a theorem about unpredictability relationships between a sequence and some of its subsequences.
We assume a class $\F$ of predictors satisfying Axioms 1-4 and a hierarchy $\F_1\subset\F_2\subset\cdots$ on this class
to be fixed. A shortened notation $I(a)=I(a;\F)$ for the unpredictability of a sequence $a$ will be used.

\begin{theorem}
\label{t:main}
We assume $a \in \{0,1\}^{\infty}$, $I(a) > 0$. For each $\gamma>0$, then either one of the five inequalities
\begin{itemize}
\item $I(P^{\nu}a) \geq I(a) + \gamma$ for $\nu = 0,1,2$,
\item $I(S^{\nu}a) \geq I(a) + \gamma$ for  $\nu = 1,2$,
\end{itemize}
holds or, for some $\tilde f \in \F$ both of the following relations
hold:
\begin{eqnarray}
\label{e:result2} \limsup_{n \to \infty} E(\tilde fa;n) \geq
\frac{I(a)}{4},\\
I(\tilde f_*a) \geq \frac{1}{2} - \frac{4\gamma}{I(a)}.
\label{e:result1}
\end{eqnarray}
\end{theorem}

Proof:
Suppose for some $a \in \{0,1\}^{\infty}$
\begin{eqnarray}
\label{e:result3}
I(P^{\nu}a) &<& I(a) + \gamma \mbox{, }\ \nu = 0,1,2\\
\label{e:result4} I(S^{\nu}a) &<& I(a) + \gamma \mbox{, }\ \nu =
1,2.
\end{eqnarray}
Then we construct a mapping $\tilde f\in\F$ such that
\eqref{e:result1} and \eqref{e:result2} hold.

For $\gamma \geq I(a)/8$, taking the constant predictor $\phi^1 \in
\F$ is sufficient for the theorem to hold. Indeed, we substitute
$I(a)/8$ into the right hand side of \eqref{e:result1} to find, $
{1}/{2} - {4\gamma}/{I(a)} \leq 0, $ but then
\begin{equation*}
I(f_*a) \geq \frac{1}{2} - \frac{4\gamma}{I(a)}
\end{equation*}
since $I \geq 0$ for all sequences. For \eqref{e:result2} we note
that $E(\phi^1a;n) = 1$ for all $n$. Since $I(a)$ is bounded above
by 1/2,  \eqref{e:result2} holds for $\tilde f=\phi^{1}$.

We fix a hierarchy of finite sets $\F_1\subset
\F_2\subset\cdots\subset\F_m\subset\cdots  $ with $\cup \F_i=\F$ and
define the notation
\begin{equation}\label{Imn}
I(a;m,n) = \min_{f \in \F_m} \frac{1}{n}\sum_{i=0}^{n-1}((f(a))_i\oplus
a_i).
\end{equation}
\begin{equation}\label{Im}
I(a;m) = \limsup_{n \to \infty} \min_{f \in \F_m}
\frac{1}{n}\sum_{i=0}^{n-1}((f(a))_i\oplus
a_i)=\limsup_{n\to\infty}I(a;m,n);
\end{equation}
hence
$$
I(a)=\lim_{m\to\infty}I(a;m)=\inf_m I(a;m).
$$
The smallest class $\F_1$ is assumed to contain predictors \eqref{psipsi} and the constant predictors $\phi^0,\phi^1$.
Suppose $0<\gamma < I(a)/8$. From assumptions \eqref{e:result3},
\eqref{e:result4}, we can fix $m_1$ such that
\begin{eqnarray}
\label{e:1.0}
I(P^{\nu}a;m_1) &<& I(a) + \gamma \mbox{, }\ \nu = 0,1,2,\\
\label{e:1.0b} I(S^{\nu}a;m_1) &<& I(a) + \gamma \mbox{, }\ \nu =
1,2.
\end{eqnarray}

It is sufficient to specify an index $m_{0}$ such that for each $m
> m_0,$ $ \alpha>0$, $n_0>0$ there is a mapping $\tilde f \in \F_{m_0}$
satisfying for some $n>n_0$
\begin{eqnarray}
\label{e:laststage1}
 E(\tilde fa;n) &>&
 \frac{I(a)}4-2\alpha, \\
\label{e:laststage2} I(\tilde f_*a;m,L) &>&
\frac{1}{2} - \frac{4\gamma}{I(a)} - \chi(\alpha),
\end{eqnarray}
where $L = nE(\tilde fa;n)$ and $\chi(\alpha)\to0$ as $\alpha\to0$.

By definition, given an $\alpha>0$, for any sequence $b$, we can choose an $N_1$ such
that $I(b;m_1,n') < I(b;m_1) + \alpha$ for all $n'>N_1$. On a finite
set $\F_{m_1}$, there must be a predictor $f$ where $E(f b\oplus
b;n')= I(b;m_1,n')$; consequently, $E(f b\oplus b;n') < I(b;m_1) +
\alpha$. Thus by \eqref{e:1.0} and \eqref{e:1.0b}, we can ensure
that if $n'$ is sufficiently large, then for some $\xi^0$, $\xi^1$,
$\xi^2$, $\eta^1$, $\eta^2 \in \F_{m_1}$:
\begin{eqnarray}
\label{e:1.2}
E(\xi^{\nu}P^{\nu}a\oplus P^{\nu}a;n') &<& I(a) + \gamma + \alpha \mbox{, } \ \nu = 0,1,2,\\
\label{e:1.21}
E(\eta^{\nu}S^{\nu}a\oplus S^{\nu}a;n') &<& I(a) + \gamma + \alpha \mbox{, } \ \nu = 1,2.
\end{eqnarray}

We construct the desired predictor $\tilde f$ using $\eta^1,
\eta^2\in \F_{m_1}$ as follows. We first use Axiom 2 to define the
predictors $c^1, c^2\in \F$ by the formulas
\begin{eqnarray}\label{x1}
P^0c^1a=P^0\phi^0a=0,\quad P^1c^1a=P^1\psi^1a=P^0a,\quad
P^2c^1a=P^2\phi^0a=0,\\
P^0c^2a=P^0\phi^0a=0,\quad P^1c^2a=P^1\psi^2a=P^2a,\quad
P^2c^2a=P^2\phi^0a=0,\label{x2}
\end{eqnarray}
where $\phi^0\in\F$ assigns the zero output sequence to any input
and the predictors $\psi^1,\psi^2\in\F$ are defined by
\eqref{psipsi}. Taking $\eta^\nu \in \F_{m_1}$ with $\nu=1,2$, by
Axiom 3 there exist $g^{\nu}_1 \in \F$ such that
\begin{equation}\label{x3}
\eta^{\nu}S^{\nu}a = P^1g_1^{\nu}a.
\end{equation}
Now we form $g_2^{\nu} \in \F$ via Axiom 2 using the predictors
$\phi^0$ and $g_1^{\nu}$:
\begin{equation}\label{x4}
P^0g_2^\nu a= P^0\phi^0a=0,\quad P^1g_2^\nu a= P^1 g_1^\nu
a=\eta^{\nu}S^{\nu}a,\quad P^2g_2^\nu a= P^2\phi^0a=0.
\end{equation}
According to Axiom 1, the predictor $g^{\nu} = c^{\nu} \oplus
g^{\nu}_2$ belongs to the class $\F$. Finally, we define the
predictor $\tilde f\in \F$ via Axiom 1 by  $\tilde f=g^1\oplus g^2$.

Remark that $g^\nu$ and $\tilde f$ belong to some sufficiently large
class $\F_{m_0}$ for any $\eta^1, \eta^2\in\F_{m_1}$. A particular
choice of $\eta^1, \eta^2$, and hence the choice of $\tilde f\in
\F_{m_0}$, depends on the value of $m$ in \eqref{e:laststage2}. In
order to specify this choice, note that Axiom 3 implies the
existence of predictors $z^1, z^2\in \F$ satisfying
\begin{equation}\label{z1z2}
P^2z^1a=\eta^1S^1a,\qquad P^2z^2a=\eta^2S^2a
\end{equation}
for any $\eta^1, \eta^2\in\F_{m_1}$. Hence, from Axiom 1 it follows
that the predictor
\begin{equation}\label{z}
z=z^1\oplus z^2\oplus \psi^2
\end{equation}
belongs to the class $\F$. Also, the predictor $f'$ defined by
\begin{equation}\label{f'def}
(f'a)_i = \left\{ \begin{array}{ll}
         (g^1a)_i, & \mbox{if $(\tilde fa)_i = 0$},\\
        (h\tilde f_*a)_{l(i)} & \mbox{if $(\tilde fa)_i = 1$}\end{array} \right.
\end{equation}
belongs to $\F$ for any $h\in\F$, according to Axiom 4.

\begin{lemma}
\label{l:axiom2'} For any $f^0, f^1, f^2 \in \F$, the predictors
$h'$ and $h''$ defined by
\begin{eqnarray}\label{star}
P^{0}h'a =  f^0 P^0  a,\quad P^{1}h'a =  P^1 f^1 a, \quad P^{2}h'a = f^2 P^2 a,\\
\label{starstar} P^{0}h''a = f^0 P^0 a, \quad P^{1}h''a = f^1 P^1
a,\quad P^{2}h''a =  P^2 f^2 a
\end{eqnarray}
belong to the class $\F$.
\end{lemma}

Indeed, Axiom 3 ensures the existence of a predictor $h^\nu\in\F$
that satisfies $ P^\nu h^\nu a = f^\nu P^\nu a $ for each $\nu=1,2$.
Now, we combine $f^0,h^1 $ and $h^2$ using Axiom 2 to obtain the
predictor $h'$ satisfying \eqref{star}. The inclusion $h''\in \F$
follows similarly. \hfill $\blacksquare$

\medskip
Given any $m$, consider a sufficiently large $m_2$ such that the
predictor \eqref{z} belongs to the class $\F_{m_2}$ for any
$\eta^1,\eta^2\in \F_{m_1}$ and the predictor \eqref{f'def} belongs
to $\F_{m_2}$ for any $h\in \F_m$, $\tilde f,g^1\in \F_{m_0}$. For
an arbitrary function $h_1 \in \F_{m_2}$, form $h_2\in \F $ from
$\xi^0$, $h_1$ and $\xi^2$ using formulas \eqref{star} of Lemma
\ref{l:axiom2'}. Consider a sufficiently large class $\F_{m_3}$ that
contains such a $h_2$ for every $h_1\in \F_{m_2}$,
$\xi^0,\xi^2\in\F_{m_1}$. From the definition of $I(a;m_3)$ it
follows that there is a sequence $n_k\to\infty$ such that
\begin{equation}
\label{e:interimIresult} I(a;m_3,n) > I(a;m_3)-\alpha\ge I(a) -
\alpha
\end{equation}
for $n=n_k, n_k+1, n_k+2$ and all $k$. Hence, there exist
arbitrarily large $n=3n'$ such that both \eqref{e:interimIresult}
holds and there are functions $\xi^\nu, \eta^\nu\in\F_{m_1}$
satisfying (\ref{e:1.2}), (\ref{e:1.21}). Consider any such $n,
\xi^\nu,\eta^\nu$ and the corresponding predictor $\tilde f\in
\F_{m_0}$ defined as described above by relations
\eqref{x1}-\eqref{x4} and $g^{\nu} = c^{\nu} \oplus g^{\nu}_2$,
$\tilde f=g^1\oplus g^2$. We will derive the desired relations
\eqref{e:laststage1}, \eqref{e:laststage2} from (\ref{e:1.2}),
(\ref{e:1.21}) and \eqref{e:interimIresult}.

Let $h^1\in\F_{m_2}$. From the relations
$$
3E(h_2a\oplus a;n) = E(P^0h_2a\oplus P^0a;n') + E(P^1h_2a\oplus P^1a;n') + E(P^2h_2a\oplus P^2a;n'),
$$
and the formulas $P^0 h_2a= \xi^0P^0 a$, $P^1 h_2a=P^1 h_1 a$, $P^2 h_2a=\xi^2 P^2 a$ defining $h_2$,
it follows that
\begin{eqnarray*}
3E(h_2a\oplus a;n) =
E(\xi^0 P^0a\oplus P^0a;n')+ E(P^1 h_1  a \oplus P^1 a;n') +
 E(\xi^2 P^2a\oplus P^2a;n').
\end{eqnarray*}
Combining this relation with \eqref{e:1.2}, we obtain
\begin{eqnarray*}
E(P^1 h_1 a\oplus P^1 a;n')  + 2(I(a) + \gamma + \alpha)  &>&  3E(h_2a \oplus a;n)
\geq 3I(a;m_3;n),
\end{eqnarray*}
where the second inequality follows since $h_2 \in \F_{m_3}$. Moreover, by \eqref{e:interimIresult}
\begin{equation*}
3I(a;m_3;n) > 3(I(a) - \alpha),
\end{equation*}
hence
\begin{eqnarray}
\label{e:1.1}
E(P^1 h_1 a\oplus P^1a;n') > I(a) -2\gamma - 5\alpha, \qquad h_1 \in \F_{m_2}.
\end{eqnarray}

 Similarly,  for each
 $h_1\in\F_{m_2}$ a predictor $h_2'$ can be formed
 by combining the predictors $\xi^0, \xi^1$ and $h_1$ according to formulas \eqref{starstar} of Lemma \ref{l:axiom2'}:
 $$
 P^{0}h'_2a = \xi^0 P^0 a, \quad P^{1}h'_2a = \xi^1 P^1 a,\quad P^{2}h'_2a =  P^2 h_1 a.
 $$
 Assuming without loss of generality that the class $\F_{m_3}$ is large enough to include
 $h_2'$ for every $h_1\in\F_{m_2}$,
 we can repeat the above argument to obtain
\begin{eqnarray}
\label{e:1.1.1}
E(P^2h_1 a\oplus P^2a;n')> I(a) -2\gamma - 5\alpha,\qquad h_1 \in \F_{m_2}.
\end{eqnarray}

Now recall the definition of $\tilde f$. It implies
\[
E(\tilde fa;n) =  \frac{1}{n} \sum_{j=0}^{n'-1} (P^{0}a)_{j}\oplus (\eta^{1}S^{1}a)_j \oplus (P^{2}a)_{j}\oplus (\eta^{2}S^{2}a)_j
\]
with $n'=n/3$. Equivalently,
\begin{eqnarray}
E(\tilde f a;n)
\label{e:combinedpredictor}
&=& \frac{1}{3}E([\eta^1S^1a\oplus \eta^2S^2a\oplus P^0a] \oplus P^2a;n').
\end{eqnarray}
Combining relations \eqref{z1z2} with the equality $P^2\psi^2 a=P^0a$, which follows from the definition \eqref{psipsi}
of $\psi^2$, we see that
$$
\eta^1S^1a\oplus \eta^2S^2a\oplus P^0a=P^2z,
$$
where $z\in \F_{m_2}$ is defined by \eqref{z}.
Hence, \eqref{e:combinedpredictor} can be rewritten as
$$
E(\tilde f a;n) = \frac{1}{3}E(P^2 z a \oplus P^2a;n')
$$
and from \eqref{e:1.1.1} it follows that
$$
E(\tilde fa;n) >     \frac{1}{3}(I(a) - 2 \gamma - 5\alpha).
$$
Together with the estimate $\gamma < {I(a)}/{8}$ this implies the desired relation \eqref{e:laststage1}.

For the second inequality, \eqref{e:laststage2},
we note that by definition of $g^\nu$, $S^\nu$,
$$
P^1g^1a=\eta^1S^1a\oplus P^0a,\quad P^1g^2a=\eta^2S^2a\oplus P^2a
$$
and $S^1a=P^0a\oplus P^1a$, $S^2a=P^1a\oplus P^2a$.
Hence, expanding the left hand side of \eqref{e:1.21}, we obtain
$$
E(\eta^{1}S^{1}a\oplus S^{1}a;n')\! =\! \frac{1}{n'}\sum\limits_{k=0}^{n'-1}(\eta^{1}S^1 a)_k \oplus (P^0 a)_k \oplus (P^1 a)_k
\!=\! \frac{1}{n'}\sum\limits_{k=0}^{n'-1}(P^1g^{1}a)_{k} \oplus (P^1a)_k
$$
and similarly
$$
E(\eta^{2}S^{2}a\oplus S^{2}a;n') = \frac{1}{n'}\sum\limits_{k=0}^{n'-1}(P^1g^{2}a)_{k} \oplus (P^1a)_k.
$$
We sum these two equations together and combine
with \eqref{e:1.21} to get
\begin{eqnarray}
\label{e:startpart2}
\frac{1}{n'}\sum_{k=0}^{n'-1}\left((P^1g^{1}a)_{k} \oplus (P^1a)_k + (P^1g^{2}a)_{k} \oplus (P^1a)_k\right) &<& 2(I(a) + \gamma + \alpha).
\end{eqnarray}
Consider the set $\J$ of indices $j < n'$ where $(P^1g^1a)_j =
(P^1g^2a)_j$ and the set $\J_c$ of indices $j< n'$ where
$(P^1g^1a)_j \ne (P^1g^2a)_j$ . From the relations
\begin{equation*}
\begin{array}{ll}
(P^1g^1a)_j\oplus (P^1a)_j + (P^1g^2a)_j\oplus (P^1a)_j =2(P^1g^1a\oplus P^1a)_j,& \qquad j\in\J,\\
(P^1g^1a)_j\oplus (P^1a)_j + (P^1g^2a)_j\oplus (P^1a)_j = 1,&\qquad
j\in \J_c
\end{array}
\end{equation*}
and \eqref{e:startpart2}, it follows that
\begin{eqnarray*}
\frac{1}{n'}\Big(\sum_{j \in \J}2(P^1g^1a\oplus P^1a)_j + \sum_{j
\in \J_c} 1\Big) &<& 2(I(a) + \gamma + \alpha).
\end{eqnarray*}
Moreover, $(P^1g^1a)_j = (P^1g^2a)_j$ is equivalent to $(P^1\tilde fa)_j = 0$, and the relation $(P^1g^1a)_j \neq (P^1g^2a)_j$ is equivalent
 to $(P^1\tilde fa)_j = 1$. Hence,
\begin{equation}\label{Ldef}
\sum_{j\in \J_c} 1 = n'E(P^1\tilde fa;n')
=nE(\tilde fa;n)
=:L,
\end{equation}
where we use the relations $P^0\tilde fa=P^2\tilde fa=0$, wich
follow from the definition of $\tilde f$. Therefore
\eqref{e:startpart2} is equivalent to
\begin{equation}
\label{e:10}
\frac{1}{n'}\Big(\sum_{j \in \J}(P^1g^1a \oplus P^1a)_j + \frac{L}{2} \Big) < I(a) + \gamma + \alpha.
\end{equation}

Let us extend the definition $(P^1\tilde fa)_j = (\tilde fa)_{3j+2}
= 1 \iff j \in \J_c$ of  the set $\J_c$ to indices $i=3j, 3j+1$. To
do this, consider the set $\J'_c$ of indices $i$ defined by
$$
\J'_c =\{i < n=3n' : (\tilde fa)_i = 1\}.
$$
Since $P^0\tilde fa=P^2\tilde fa=0$, we see that $i\in
\J_c'$ if and only if $i=3j+2$ with $j\in \J$, hence for any sequence $b$
\begin{equation}\label{b*}
\sum_{j\in \J_c} (P^1b)_j = \sum_{i\in \J_c'} b_i.
\end{equation}
Now, recall that for any $h \in \F_{m}$, using Axiom 4, we can
construct the function $f'\in\F_{m_2}$ defined by \eqref{f'def}.
Applying the identity \eqref{b*} to the sequence $b=f'a\oplus a$, we
obtain,
$$
\sum_{j\in \J_c} (P^1f'a\oplus P^1a)_j = \sum_{i\in \J'_c} (f'a\oplus a)_i
= \sum_{i \in \J'_c} (h\tilde f_*a)_{l(i)}\oplus a_i,
$$
where the second equality follows from the definition of $f'$ and
$\J'_c$. (The notation $l(i)$ is introduced in Axiom 4; $l(i)$ is
the number of 1's in the sequence $\tilde fa$ up to, but not
including, the digit $(\tilde fa)_i$.)\ As $\tilde f_*a$ is, by
definition, the subsequence selected from $a$ whenever $(\tilde
fa)_i = 1$,
$$
a_i=(\tilde f_* a)_{l(i)},\qquad i\in \J_c',
$$
hence
\begin{equation}\label{zzz}
\sum_{j\in \J_c} (P^1f'a\oplus P^1a)_j = \sum_{i \in \J'_c} (h\tilde f_*a)_{l(i)}\oplus (\tilde f_* a)_{l(i)}=
\sum_{k=0}^{L-1} (h\tilde f_*a)_{k}\oplus (\tilde f_* a)_{k}.
\end{equation}
Here $L$ is the cardinality of the set $\J_c'$, which is equal to the cardinality of the set $\J_c$, hence $L$ is
defined by formulas \eqref{Ldef}.
Now note that if $j\in \J$, then $(P^1\tilde f a)_j=(\tilde f a)_{3j+2}=0$,
hence $(f' a)_{3j+2}=(g^1 a)_{3j+2}$, that is $(P^1 f'a)_j=(P^1 g^1 a)_j$ for $j\in\J$.
Therefore
\begin{equation}\label{zzz'}
\sum_{j\in \J} (P^1f'a\oplus P^1a)_j = \sum_{i \in \J} (P^1 g^1 a\oplus P^1a)_j.
\end{equation}
Summing \eqref{zzz} and \eqref{zzz'}, we obtain
$$
E(P^1 f'a\oplus P^1a; n') =\frac{1}{n'}\Big( \sum_{k=0}^{L-1} (h\tilde f_*a)_{k}\oplus (\tilde f_* a)_{k}+\sum_{j \in \J} (P^1 g^1a)_j \oplus (P^1a)_j \Big),
$$
hence \eqref{e:1.1} implies
\begin{equation}
\label{e:11}
\frac{1}{n'}\Big( \sum_{k=0}^{L-1} (h\tilde f_*a)_{k}\oplus (\tilde f_* a)_{k}+\sum_{j \in \J} (P^1 g^1a)_j \oplus (P^1a)_j \Big) > I(a) -2\gamma - 5\alpha.
\end{equation}
Furthermore, subtracting \eqref{e:10} from \eqref{e:11} we arrive at
\begin{equation*}
\frac{1}{n'} \Big(\sum_{k = 1}^{L-1} (h\tilde f_*a)_k \oplus (\tilde f_*a)_k - \frac{L}{2} \Big) > -3\gamma - 6\alpha.
\end{equation*}
Equivalently,
$$
\frac{1}{L} \sum_{k = 1}^{L-1} (h\tilde f_*a)_k \oplus (\tilde f_*a)_k> \frac12 -\frac{n'(3\gamma+6\alpha)}L
=\frac12 -\frac{n(\gamma+2\alpha)}L.
$$
These relations combined with \eqref{e:laststage1} and \eqref{Ldef} imply
\begin{equation}\label{kkk}
\frac{1}{L} \sum_{k = 1}^{L-1} (h\tilde f_*a)_k \oplus (\tilde f_*a)_k> \frac12 -
\frac{\gamma+2\alpha}{\frac{I(a)}4-2\alpha}=
\frac12 -
\frac{4\gamma}{I(a)} - \chi(\alpha)
\end{equation}
with $\chi(\alpha)\to0$ as $\alpha\to0$. Finally, as \eqref{kkk}
holds for an arbitrary $h \in \F_{m}$, we infer the estimate
\eqref{e:laststage2}. This completes
 the proof of the theorem.\hfill
$\blacksquare$

\subsection{Independence}
We combine the above theorem with an idea of independence, which has a certain analogy to the idea of independence in probability theory.
\begin{definition}
We say that a sequence $a$ consists of $\F$-independent quantities (or, shortly, that $a$ is ${\mathcal F}$-independent)
if, for any $f \in \F$,
\begin{equation*}
I(f_* a) = I(a).
\end{equation*}
\end{definition}

{\bf D: Alexei, would we need more discussion of F-independence at this point?}

$\F$-independence enables the following theorem.

\begin{theorem}\label{t:independence}
Suppose a sequence $a$ consists of $\F$-independent quantities.
Define the sequence $b^{\nu}$ with $\nu=1,2$ by
$b^{\nu}_i = a_{3i+\nu-1} \oplus a_{3i +2}$ for $i\ge1$. Then
the following inequality holds for at least one $b^{\nu}$
\begin{equation}
\label{c:independence}
I(b^{\nu}) \geq I(a) \left(1 + \frac{1 - 2I(a)}{5}\right).
\end{equation}
Hence, $I(b^{\nu})>I(a)$ for at least one $b^{\nu}$ whenever $0<I(a)<1/2$.
\end{theorem}

Proof: Relation \eqref{c:independence} is trivial for $I(a)=0$,
hence assume $I(a)>0$. We first prove that $I(a) = I(P^{\nu}a)$ for
$\F$-independent sequences. We choose the predictor $f =
001001\ldots$. This can be formed from the constant predictors
$\phi^0$ and $\phi^1$ and use of Axiom 2, thus $f \in \F$. Then
since $a$ is $\F$-independent
\begin{equation*}
I(a) = I(f_*a) = I(P^{0}a).
\end{equation*}
Similar constructions for $f$ provide the result for other values of
$\nu$. We note that  $b^{\nu} = S^{\nu}a$. We now apply Theorem
\ref{t:main} with $\gamma = I(a)\left(\frac{1 - 2I(a)}{5}\right)$.
Since $I(a) = I(P^{\nu}a)$, the relations $I(P^{\nu}a)\ge I(a) +
\gamma$, can not hold. Thus either, for at least one $b^{\nu}$ we
have
\begin{equation}
\label{e:stratalast} I(b^\nu) = I(S^{\nu}a) \ge I(a)\left(1 +
\frac{1 - 2I(a)}{5}\right)
\end{equation}
or inequalities \eqref{e:result1}, \eqref{e:result2} hold for some
$\tilde f$. In the latter case,
\begin{equation*}
I(a) = I(\tilde f_*a) \geq \frac{1}{2} - \frac{4\gamma}{I(a)},
\end{equation*}
since $a$ is $\F$-independent, and substituting in $\gamma$ gives
$$
I(a) \geq \frac{1}{2} - \frac{4 I(a) \left(\frac{1 - 2I(a)}{5}\right)}{I(a)}
\geq \frac{1}{2} - 4\left(\frac{1 - 2I(a)}{5}\right).
$$
This implies $1/2\ge I(a)$,
which is a contradiction if $I(a) \neq \frac{1}{2}$. Thus \eqref{e:stratalast} holds, and the theorem is proved.
\hfill
$\blacksquare$

{\bf D: Alexei, the above proof does not work for $I(a)=1/2$, otherwise OK}.

{\bf F: I had a think about this and couldn't think of an obvious way to make it work. Am I missing a trivial argument that $I(a)=1/2$?}.

We can compare this result to results in the classical probability
formalism. Suppose we have a sequence of independent identically
distributed random variables $X_i$ taking binary values 0 with
probability $p$ and 1 with probability $q= 1-p$.  Now for individual
realisations of such sequences, we show that almost all (in the
probabilistic sense) will have unpredictability $I(a) = \min\{p,q\}$
which is achieved by one of the constant predictors $\phi^1$ or
$\phi^0$.

\begin{theorem} Consider the set of sequences generated by realisations of a sequence of independent identically distributed binary random variables
$X_i$ with $\P[X = 0] = p$, and $\P[X=1]=q$ for $X=X_i$. Almost
every realisation, $x$ has an unpredictability value $I(x) =
\min\{p,q\}.$
\end{theorem}

Proof:
We note first that an upper bound on $I(x)$ is achieved by one of the constant functions $\phi^0,\phi^1$. By the strong law of large numbers,
$$
\lim_{n \to \infty} \frac{1}{n}\sum_{i=0}^{n-1} X_i = \E[X] = q
$$
almost surely. Similarly,
$$
\lim_{n \to \infty} \frac{1}{n}\sum_{i=0}^{n-1} X_i \oplus 1 =
\E[X\oplus 1] = \sum_{x=0,1}(x \oplus 1) \P[X=x] = p
$$
almost surely. Hence
\begin{eqnarray*}
I(x) \leq
\min \{p,q\}
\end{eqnarray*}
for almost every realisation $x$. For the lower bound, consider
\begin{eqnarray*}
\P[(f(X))_i\oplus X_i = 1] & = & \P[X_i=0]\P[(f(X))_i = 1] + \P[X_i=1]\P[(f(X))_i = 0]\\
&=& p\P[(f(X))_i = 1] + q(1 - \P[(f(X))_i=1])\\
&=& (p-q)\P[(f(X))_i = 1] + q,
\end{eqnarray*}
where we use the fact that the events $X_i=0$ and $(f(X))_i=1$ are
independent, as the events $X_i=1$ and $(f(X))_i=0$ are, because
$(f(X))_i$ is a function of the variables $X_1,\ldots, X_{i-1}$ only
and hence $X_i$ are $(f(X))_i$ are independent. Similarly,
\begin{equation*}
\P[(f(X))_i\oplus X_i = 1] = (q-p)\P[(f(X))_i = 0]+p,
\end{equation*}
and thus for each predictor $f$
\begin{equation}
\label{e:probpq}
\P[(f(X))_i\oplus X_i = 1] \geq \min\{p,q\}.
\end{equation}
{\bf D: I did not get the rest of the proof from this point.}

Now, we can write:
\begin{eqnarray*}
\E\Big[(f(X))_i\oplus X_i\Big] = \P[(f(X))_i\oplus X_i = 1]
\end{eqnarray*}
Thus by \ref{e:probpq}, and by the strong law,
\begin{eqnarray*}
\min\{p,q\} \leq \E\Big[(f(X))_i\oplus X_i\Big] = \lim_{n \to \infty} \frac{1}{n} \sum_{i=0}^{n-1} f(a)_i \oplus a_i
\end{eqnarray*}
on a set of sequences of measure 1. But this is true for all $f$, so we can write
\begin{eqnarray*}
\min\{p,q\} \leq \lim_{n \to \infty} \inf_{f \in \
F} \frac{1}{n} \sum_{i=0}^{n-1} f(a)_i \oplus a_i = I(a)
\end{eqnarray*}
which is true on a set of sequences of measure 1. Thus we have established both bounds, hence
\begin{eqnarray*}
I(a) = \min\{p,q\}
\end{eqnarray*}
on a set of sequences of measure 1.
\hfill
$\blacksquare$

If we examine the probability distribution on the sequence $b =
a_{3i}\oplus a_{3i-1}$, we find each $b_i$ takes value 0 with
probability $p^2 + (1-p)^2 = 2p^2 -2p +1$ and takes value 1 with
probability $2p(1-p) = 2p - 2p^2$. So using the constant predictors,
$\phi^0$ and $\phi^1$, by a similar argument to above, we can
guarantee
$$
I(b)= \min\{1 - (2p -2p^2), 2p - 2p^2\}
 = 2p - 2p^2
$$
since $2p - 2p^2 \le 1/2$ for all $p \in [0,1]$. Now if $p <1/2$, $I(a) = p$ and $I(b) = 2I(a) - 2I(a)^2$. If $p>1/2$, $I(a) = 1-p$ and $I(b) = 2p - 2p^2 = 2(1-p) - 2(1-p)^2 = 2I(a) - 2I(a)^2$. Thus we can write this relation in the form of \eqref{c:independence}, i.e.,
$$
I(b) = 2I(a) - 2I(a)^2
= I(a)(1 + (1-2I(a)).
$$
This is a more exact result than \eqref{c:independence}, though obtained from more restrictive conditions.
It implies that for almost every Bernoulli sequence $a$ with $0<p<1/2$, $q=1-p$
$$
I(b)>I(a),
$$
i.e., the simple operation producing the sequence $b_i=a_{3i}\oplus a_{3i-1}$ increases the unpredictability.
The authors do not know whether the bound \eqref{c:independence} obtained in Theorem \ref{t:independence}
through the condition of $\F$-independence is tight.

\section{Appendix 1: Proof of Theorem \ref{tt2}}
We first show how to construct a sequence $a$ with $I(a)=1/2$.
Consider a particular predictor $f_1 \in \F$, acting on a finite
sequence of length $n$. If $I(a;f_1,n) = 0$, then
\begin{equation*}
(f_1(a))_i= a_i
\end{equation*}
for all $i = 1,\ldots,n$. That is, the sequence $f_1(a)$ is
completely defined - there is only one sequence with $I(a;f_1,n) =
0$. For $I(a;f_1,n) = 1/n$, then $(f_1(a))_i\neq a_i$ occurs at one
and only one element of $a$. Thus there are $n$ sequences with
$I(a;f_1,n) = 1/n$. In general for $I(a;f_1,n) = k/n$,
$(f_1(a))_i\neq a_i$ can occur in ${n \choose k}$ combinations,
hence $f_1$ predicts ${n \choose k}$ sequences with $I(a;f_1,n) =
k/n$.

We now consider, for large $n$, the class of sequences,
$\#A_{f_1,n,\epsilon}$, with
\begin{equation}\label{anp}
|I(a;m,n)-1/2|<\varepsilon.
\end{equation}
The cardinality of this class is
\begin{equation*}
\label{d:typicalset} \#A_{f_1,n,\epsilon} = \sum_{k = \lceil n/2 -
n\epsilon \rceil}^{k = \lfloor n/2 + n\epsilon \rfloor} {n \choose
k}.
\end{equation*}

The following lemma is a variation on the De Moivre - Laplace theorem, see for example \cite{U37}, see also the original version by De-Moivre in \cite{D00}.

\begin{lemma} For any $\epsilon, \delta > 0$ there is an $N_1=N_1(\delta)$ such that for all $n\ge N_1$
\begin{equation*}
\#A_{f_1,n,\epsilon}>2^{n - \delta}.
\end{equation*}
\end{lemma}

For any finite set of predictors, $\F_m=\{f_1,\ldots,f_p\}$, the set
of sequences with unpredictability satisfying \eqref{anp} is
\begin{equation*}
\bigcap_{i} A_{f_i,n,\epsilon}
\end{equation*}
which has cardinality
\begin{equation}
\# \bigcap_{i} A_{f_i,n,\epsilon}> 2^n - \sum_{i=1}^{p}(2^n -
2^{n-\delta})>2^n\bigl(1-p(1-e^{-\delta})\bigr)
\label{e:almostfull}
\end{equation}
for all $n\ge N=N(\delta)$, where $N =\max(N_1,\ldots,N_p)$. For a
sufficiently small $\delta$, we see that the set of sequences with
$|I(a;m,N) -\frac{1}{2}|< \epsilon$ is non-empty for $n\ge N$ - in
fact, it is almost the full set (not unlike the ``typical set'' in
the information theory sense).

Let $a'$ be an arbitrary block of length $|a'|$. There are
$2^{n-|a'|}$ sequences of length $n>|a'|$ beginning with $a'$. Now,
given $a'$, $\F_m$ and $\epsilon>0$, if $\delta$ is sufficiently
small, then for any $n$
\begin{equation*}
2^{n-|a'|} + 2^n\bigl(1-p(1-e^{-\delta})\bigr)> 2^n
\end{equation*}
and hence \eqref{e:almostfull} implies that there exist sequences
$a$ beginning with block $a'$ for which we can choose an
$N=N(\epsilon, |a'|)$ such that $|I(a;m,N) - \frac{1}{2}|<
\epsilon$. Consequently, we can choose blocks $a^1, a^2, \ldots$,
with lengths $N_1, N_2 - N_1, N_3 - N_2 \ldots$ respectively, and
guarantee that these blocks satisfy
$$
|I(a^1a^2\ldots a^m;m,N_m) - {1}/{2}|< \epsilon_m
$$
with $\epsilon_i = 2^{-i}\epsilon_1$ for all $m\ge 1$.

For all $j$, $n$ and $a$, the inclusion $\F_{j-1}\subset\F_j$
implies $I(a;j-1,n) \geq I(a;j,n)$. Thus for a given class $\F_m$,
at sequence lengths $N_1,N_2,\ldots, N_j$
$$
\limsup_{j \to \infty} I(a^1a^2\ldots a^j;m,N_j) \geq \limsup_{j \to
\infty}I(a^1a^2\ldots a^j;j,N_j) \geq  \lim_{j \to \infty}
\Big(\frac{1}{2} - \epsilon_j\Big) =\frac{1}{2}.
$$
Define $a=a^1a^2a\ldots$. We know that at points $n = N_j$,
\begin{equation*}
I(a;m,n) = I(a^1a^2\ldots a^j;m,N_j)
\end{equation*}
and
\begin{equation*}
I(a;m) = \lim_{n \to \infty} \sup I(a;m,n)\geq \lim_{j \to \infty}I(a^1a^2\ldots a^j;m,N_j) \geq \frac{1}{2}.
\end{equation*}
Since $\phi_0, \phi_1 \in \F$, $I(a;m,n)$ is also bounded above by
$1/2$ and hence $
 I(a;m) = \frac{1}{2}$
for all sufficiently large  $m$. Consequently,
\begin{equation}
I(a)= \lim_{m \to \infty} I(a,m) = \frac{1}{2}. \label{e:halfexists}
\end{equation}

Now we show how to construct a sequence with any unpredictability
$I_0<\frac{1}{2}$.

We first extract a slightly stronger statement from the preceding arguments;
for an unspecified predictor class of given cardinality, we require that we can generate a sequence of
high unpredictability within a guaranteed number of digits. Specifically, the next lemma
follows directly from \eqref{e:almostfull}.
\begin{lemma}
\label{l:extracted} For any $p$ and $\epsilon >0$, there exists an
$N$ such that for each set $\tilde \F$ of predictors of size
$\#\tilde\F\le p$, there exists a finite sequence $a$ of length $N$
such that
\begin{equation*}
I(a;\tilde\F,N) > \frac{1}{2} - \epsilon.
\end{equation*}
\end{lemma}

This allows us to prove the following statement.
\begin{lemma}
\label{l:blocks} For any predictor class $\F_m$, any $\epsilon>0$
and any finite sequence $a$ of length $n$, there exist an $N'$ and
blocks $b$ of any length $N>N'$ such that when block $b$ is
concatenated with sequence $a$,
$$
\inf_{f\in \F_m} \frac{1}{N}\sum_{i=0}^{N-1}f(ab)_{i+n}\oplus b_i >
\frac{1}{2} - \epsilon.
$$
Moreover, $N'$ is independent of the length $n$ of $a$.
\end{lemma}

{\bf D: 1. Lemmas 6.2 and 6.3 look very similar; maybe the first
follows from the second one.

F: We actually use the first one to prove the second one.

2. Moreover, both of these lemmas look very similar to the
statements and argument used at the beginning of the proof on page
16.

F:We use the arguments on page 16 to prove the lemmas. page 16$->$ lemma 6.2 $->$6.3.

3. There is no reference to Lemma 6.2 further. There is no reference
to Lemma 6.3 in this appendix either --- the first reference appears
in Appendix 3.

F: It's used directly after, I've put in the references explicitly.

4. Hence, can we formulate just one lemma at the beginning of this
proof and refer to it systematically? The structure, as it is, seems
somewhat confusing to me.

5. I did not work through the rest of the proof, i.e. proving the
unpredictability values between 0 and 1/2, feeling that this
structural thing should be sorted out first.}

Proof: We can consider a finite sequence $a$ of length $n$ as a mapping on the space of predictors, $a: \F \to \F$, in the following manner:
$$
a(f(b))_i= f(ab)_{i+n}
$$
for $b \in \{0,1\}^{\infty}$. Let $a(\F_m)$ denote the set of
predictors obtained by $a$ acting on each predictor in $\F_m$. Now,
for any $\epsilon >0$ one can find an $N'$ such that for each $N>N'$
there exists a sequence $b$ of length $N$ such that
\begin{equation}
\inf_{f\in\F_m}\frac{1}{N}\sum_{i=0}^{N-1}f(ab)_{i+n}\oplus b_i =
\inf_{g \in a(\F_m)} \frac{1}{N}\sum_{i=0}^{N-1}g(b)_i \oplus b_i
> \frac{1}{2} - \epsilon \label{e:constructahalf}
\end{equation}
which follows from the same arguments leading to
\eqref{e:halfexists}. Hence, there exist sequences of
unpredictability $1/2$ for any set of predictors $\F$. Independence
of $N$ from $n$ (the length of $a$), follows from the fact that
$\#a(\F_m) \leq \#\F_m$, and Lemma \ref{l:extracted}. \hfill
$\blacksquare$

We now use lemma \ref{l:blocks} to demonstrate existence of sequences with arbitrarily chosen unpredictability value. Consider the change in $I(a)$ if we add a block $b^1$ obtained from lemma \ref{l:blocks}:
$$
I(a) - I(ab^1) = I(a) - \frac{nI(a) + (\frac{1}{2} - \epsilon)(N)}{n+N}.
$$
This tends to zero as $n$ tends to $\infty$. Specifically, for any
arbitrary $\delta > 0$ we can find an $n'$ such that for all $n>n'$
adding a block $b^1$ will result in a change of less than $\delta$.
If we take a sequence of $k$ zeroes, $a=000\ldots$, and form the
infinite sequence
$$
a' = ab^1b^2b^3\ldots,
$$
then $I(a';m) = \frac{1}{2} - \epsilon$.  $I(a;m;n)$ starts at zero,
and we choose $k$ large enough such that we increase I in steps of
size less than $\delta/m$. At some point
$$
I(ab^1\ldots b^r) < I_0 < I(ab^1\ldots b^rb^{r+1})
$$
and the sequence truncated at block $b_r$ has
$$
I_0 - \frac{\delta}{m} < I(ab^1 \ldots b^r) < I_0.
$$
Now we construct a sequence $c$ with $I(c) = I_0$. First construct a
block $c^1$ using the previous construction for m=1. Then choose a
block, $a^2$, of zeros such that we are within $\epsilon$ of zero
(and choose $\epsilon<I_0$), and long enough that the block size of
the above construction with $m=2$ will be less than $\delta/m$. We
then construct $c^2$ by the above method but with $m=2$. Continuing
this process we generate the sequence $c=a^1c^1a^2c^2a^3c^3\ldots$
$$
I_0 - \frac{\delta}{m}<I(a^1c^1a^2c^2\ldots a^m c^m;m)<I_0.
$$
We now show a lower bound on I(c). For any fixed $m$ we can find $n = |a^1c^1a^2c^2\ldots a^m c^m|$ such that
$$
I(c;m;n) > I_0 - \frac{\delta}{m}
$$
Also, at the end of each block $b^i$ in $c$ with $n'>n$,
$$
I(c;m;n') > I(c;m+j;n') > I_0 - \frac{\delta}{m+j}
$$
for all $j>0$, up to where $|a^1c^1\ldots c^{m+j}| = n'$. Thus
$$
\lim_{n \to \infty}\sup I(c;m;n) \geq \lim_{n \to \infty} I_0 - \frac{\delta}{m+j} = I_0
$$
Now the upper bound on $I(c)$. We examine $I(c;m;n)$ at an arbitrary $c^i$ block, with $i\geq m$ We know the value of
unpredictability truncated at subblocks $b^j$ within $c^i$ is increasing in steps of $\delta/m$. Thus the
highest unpredictability occurs in the last $b^j$ block. The increase in $I$ from the beginning of $b^j$ to the end is bounded by $2\delta /i$.
But the value at the end, $I(a^1c^1\ldots c^i;m) < I_0$, thus the value of $I(c;m)$ over $c^i$ is bounded by $I_0 + 2\delta/i$.

Now consider the start of the $c^{i+1}$ block. Suppose the following case: that the zero predictor, $\phi^0$ has value $I_0 + 2\delta/i$.
Then as we examine the unpredictability at increasing digits of $c^{i+1}$ the unpredictability increases at most to $I_0 + \delta/2i$
(the case where the best predictor predicts continuously wrong, until crossing with the $\phi^0$ predictor which is
predicting continuously correct within $c^{i+1}$). In general for any value of $I \in [I_0 - \delta/i, I_0 + 2\delta/i]$,
the value of the increase is bound by the decreasing value of the $\phi^0$ predictor, which is bounded by a monotonic decrease from
$I_0 + 2\delta/i$. Thus
$$
\lim_{n \to \infty}\sup I(c;m;n) \leq \lim_{n \to \infty} I_0 + 2\delta/i = I_0
$$
since $i \to \infty$ as $n \to \infty$. This holds for all $m$, and hence $I(c) = I_0$.
\hfill
$\blacksquare$

\section{Appendix 2: Examples of predictor classes}
Here we show that two classes of predictors, the finite state machines and the Turing machines,
satisfy the set of Axioms 1--4 stated in Section 3. Hence, the measure of unpredictability
defined by each of these classes satisfies the conditions of Theorems \ref{t:main}, \ref{t:independence}.

\subsection{Finite state machines} There are a number of
alternative definitions of a finite state machine. The idea of a
finite state machine has roots in both computer science and
linguistics, in particular an area known as formal language theory.
Originally investigated in the 60's, they have more recently found
use as a method of representation of the control logic and program
flow in software design. They are less well known for their
interpretation as predictors, which is what we will use them for.
When we refer to a finite state machine, we mean the
definition of a Moore machine.

\begin{definition}
A Moore machine is a sextuple,
$$
M  = (X,Y,S,s_0,\lambda,\delta)
$$
where
\begin{itemize}
\item $X$ is a finite set, the set of inputs (here restricted to $\{0,1\}$),
\item $Y$ is a finite set, the set of outputs (here restricted to $\{0,1\}$),
\item $S$ is a finite set, the set of states,
\item $s_0$ is a an element from S - the initial active state of the machine,
\item $\lambda: S\times X  \to S$, is the state transition function,
\item $\delta: S \to Y$, is the output function.
\end{itemize}
\label{d:automata}
\end{definition}
We will simplify our working conditions in this study by always
working with binary machines, that is both $X$ and $Y$ are
$\{0,1\}$.

If we input any sequence to a finite state machine, the output
sequence,
$$
\delta(s_0),
\delta(\lambda(s_0,a_0)),\delta(\lambda(\lambda(s_0,a_0),a_1)),\ldots
$$
defines a function on both $\{0,1\}^*$ (all finite binary sequences)
and $\{0,1\}^{\infty}$. We can consider this sequence as predictions
of the sequence $a_i$ with the property of causality - $\delta(s_0)$
is our prediction for $a_0$, $\delta(\lambda(s_0,a_0))$ is our
prediction for $a_1$ and so on. Thus a finite state machine can be
considered as a predictor.

We note that a natural hierarchy exists for finite state machines
--- they can be ordered by the number of states they contain.

\begin{theorem}
The class of all finite state machines
satisfies Axioms 1 -- 4.
\end{theorem}

Proof:

{\bf Axiom 1 (Summation).} Given finite state machines $f^0,f^1$
with $Q^0$ and $Q^1$ states, respectively, we construct the machine
$f=f^0\oplus f^1$ as follows. Define $Q^0Q^1$ states of $f$. We
associate each state in $f$ with a state in $f^0$ and a state in
$f^1$. Accordingly, we label the states in $f$ by the pair $s^0_i
s^1_j$. Suppose $\lambda^0,\lambda^1$ are the transition functions
for machines $f^0$ and $f^1$, and suppose $\delta^0,\delta^1$ are
the output functions for machines $f^0,f^1$. We define the
transitions of $f$ as
$$
\lambda(s^0_i s^1_j,a_k)= \lambda^0(s^0_i,a_k)\lambda^1(s^1_j,a_k),
$$
and define the output as
$$
\delta(s^0_is^1_j)= \delta^0(s^0_i) \oplus \delta^1(s^1_j).
$$
This machine with the initial state $s_0^0s_0^1$ behaves as the
desired predictor with $Q^0Q^1$ states.

\medskip {\bf Axiom 2 (Interleaving). }
Consider the state machines $f^0,f^1,f^2$, with $Q^0,Q^1,Q^2$ states
respectively. Form a new machine $f$  with $3Q^0Q^1Q^2$ states,
labelling each state by $\gamma s^0s^1s^2$, where $\gamma$ takes
values $0, 1$ or $2$ and $s^i$ is a state of the machine $f^i$.
Define the state transition and output functions of $f$ by
\begin{eqnarray*}
\lambda(0s^0s^1s^2,a_k) &=& 2\lambda^0(s^0,a_k)\lambda^1(s^1,a_k)\lambda^2(s^2,a_k)\\
\lambda(2s^0s^1s^2,a_k) &=& 1\lambda^0(s^0,a_k)\lambda^1(s^1,a_k)\lambda^2(s^2,a_k)\\
\lambda(1s^0s^1s^2,a_k)&=&0\lambda^0(s^0,a_k)\lambda^1(s^1,a_k)\lambda^2(s^2,a_k)
\end{eqnarray*}
and
$$
\delta(0s^0s^1s^2) = \delta^0(s^0),\quad \delta(1s^0s^1s^2) =
\delta^1(s^1),\quad \delta(2s^0s^1s^2) = \delta^2(s^2),
$$
where $\lambda^i$, $\delta^i$ are the state transition function and the output function of the machine
$f^i$.
This machine with the initial state $0s_0^0s_0^1s_0^2$ behaves as
$f$ constructed via Axiom 2.

\medskip
{\bf Axiom 3 (Subsequences). } We first construct the machine $h^0$
satisfying $P^0h^0a=fP^0a$ as required in Axiom 3. This is
accomplished by inserting two extra dummy states for each state in
$f$. More precisely, for every state $s$ in $f$, we define the
states $0s, 1s, 2s$ in $h^0$. Define the transition function for
$h^0$ as
$$
\lambda'(2s,a_k) = 0\lambda(s,a_k),\quad \lambda'(0s,a_k) = 1s,\quad
\lambda'(1s,a_k) = 2s
$$
and the output function as
$$
\delta'(2s)= \delta(s)
$$
with output for $0s,1s$ defined arbitrarily; here $\lambda$ and
$\delta$ are the transition and output function for $f$. Define the
starting state in $h^0$ as $2s_0$ where $s_0$ is the starting state $f$.
This completes the construction of $h^0$. Machines $h^1$ and $h^2$ can be constructed
in a similar manner.

Now we construct the machine $f^1$ satisfying $P^0f^1a=fS^1a$ by inserting an extra state at each $0s$ position.
We thus require four states $0s, 1s, 2s, 3s$ in the machine $f^1$ for each state $s$ in $f$.
The state transition function $\lambda'$ and the output function $\delta'$ of $f^1$ are defined by
\begin{align*}
&\lambda'(0s,a_k) = 1s,\ \ \
\lambda'(1s,0)= 2s,\ \ \
\lambda'(1s,1) =3s,\\
&\lambda'(2s,0) = \lambda'(3s,1)=0\lambda(s,0),\ \ \
\lambda'(2s,1) = \lambda'(3s,0)=0\lambda(s,1)
\end{align*}
and
$$
\delta'(0s)=\delta(s)
$$
with $\delta'$ arbitrarily defined on the states $1s, 2s, 3s$. The starting state of $f^1$ is $0s_0$.

If $f$
has $Q_f$ states, this machine satisfies the desired constraint with
$4Q_f$ states. Constructing machines $f^2,g^1,g^2$ to satisfy the other three
constraints for Axiom 3 is done in a similar fashion, each new
machine requiring $4Q_f$ states.

\medskip
{\bf Axiom 4 (Switching). } Given machines $f^0,f^1,f^2$ with
$Q^0,Q^1,Q^2$ states respectively, we define a state machine $f$
with $Q^0Q^1Q^2$ states. We label the states of $f$ by
$\gamma s^0s^1s^2$, corresponding to the sets of states $s^0,s^1,s^2$ of the
machines $f^0,f^1,f^2$, where $\gamma=0$ if $\delta^0(s^0)=0$ and
$\gamma=1$ if $\delta^0(s^0)=1$. Hence,
 the composite machine $f$ is defined by examining whether the output $\delta^0(s^0)$ of $f^0$ is zero or one. If zero, we  output according to the machine $f_1$, and update the states of the machines $f_0$ and $f_1$. If $\delta^0(s^0)=1$, then we output according to the machine $f_2$, and update the states of the machines $f_0,f_1$ and $f_2$. Thus we define the transition and the output functions of $f$ by
\begin{eqnarray*}
\lambda(0s^0s^1s^2,a_k
) &=& \lambda^0(s^0,a_k) \lambda^1(s^1,a_k) s^2,\\
\lambda(1s^0s^1s^2,a_k
) &=& \lambda^0(s^0,a_k) \lambda^1(s^1,a_k) \lambda^2(s^2,a_k),\\
\delta(0s^0s^1s^2)&=& \delta^1(s^1),\\ \delta(1s^0s^1s^2)&=& \delta^2(s^2).
\end{eqnarray*}
This machine with the initial state $\delta^0(s^0_0)s^0_0s^1_0s^2_0$ satisfies Axiom 4 by construction.
\hfill $\blacksquare$

\subsection{Turing machines}
We provide another example of a class of predictors based on Turing
machines - more specifically, the recursive predicate functions
(defined below). In another language these are the set of all
computable predictors. We first define recursive functions, which we
do via the definition of a Turing machine.

\begin{definition}
A Turing machine consists of a tape and a finite control. The tape
consists of an infinite amount of cells, $c_i$, $i\in \Z$ each of
which contains either a zero, a one, or a blank symbol. The finite
control is a finite state machine, which reads values from the tape
as input. Time, $t=0,1,2,\ldots$, is the steps of the state machine and at
time $t=0$ the state machine is positioned to read cell $c_0$ as
input. The output of the state machine is to either
\begin{itemize}
\item Move left  - If finite control is positioned at cell $c_i$, then prepare to read cell $c_{i-1}$,
\item  Move right  - If finite control is positioned at cell $c_i$, then prepare to read cell $c_{i+1}$,
\item If finite control is positioned at cell $c_i$, then rewrite the value of $c_i$ to either zero, one, or blank.
\end{itemize}
At time $t= 0$, the tape has a continuous finite sequence of zeros
and ones stretching from $c_0$ to the left, and all other cells are
blank. This is known as the input, or the program. Lastly, the
finite control has a special halting state; if this state is reached
the machine reads no more input and halts. The state of the tape
after the machine halts is the output of the Turing machine.
\end{definition}

\begin{definition}
A self delimiting version of a finite sequence $a$, denoted
$\overline{a}$ is the sequence $a$ concatenated together with a
prefix which encodes the length of $a$, $l(a)$.
\end{definition}

For example, a simple scheme for describing the length of $a$ is adding $l(a)$
1's to start of the sequence, followed by a zero to describe the
end, that is
$$
\overline{x} = 1^{l(x)}0x.
$$
Here we know the length of $a$ by counting the number of ones up to
the first zero. After that zero, we can be sure that the string $a$
is beginning. Other more efficient schemes exist.

A partial function is a function which is not necessarily defined
for all values of its domain.
We can associate a partial function with each Turing machine.

\begin{definition}
Represent the $n$-tuple of integers $(x_1, \ldots, x_n)$ by a single
binary string consisting of a concatenation of self-delimiting
versions of all the $x_i$'s. Use this as input to a Turing machine.
The integer represented by the binary string that occupies the tape
at the time of the machine halting is the value of the partial
function associated with the Turing machine, $p:\N^n \to \N$. These
functions are the \textit{partial recursive} or $computable$
functions.
\end{definition}

\begin{definition}
If the associated Turing machine halts for all inputs, the function
is known as $recursive$ function.
\end{definition}

We examine functions with a restriction of the range to $\{0,1\}$
--- these are known as predicate functions, \cite{livitanyi}. Now
predicate functions which are also recursive output a $1$ or $0$ for
all inputs of finite length, thus for each recursive predicate
function, $R$ say, we can define a predictor:
$$
(f(a))_{i+1} = R(a_1 \ldots a_i).
$$
The first digit of the prediction is arbitrary. We will call these
predictors \textit{Recursive predictors}. We will consider the
unpredictability definition with respect to the set of all recursive
predictors.

{\bf D: Finn, Alexei, I didn't quite get the definition of the predictor. Definitions 7.5, 7.6 define a function
$p: N^n\to N$. How a function $f: \{0,1\}^\infty\to \{0,1\}^\infty$ is defined based on $p$?
Why $f$ is causal? }

\begin{theorem}
The set of all recursive predictors is closed under Axioms 1 -- 4.
\end{theorem}

We sketch the proof, omitting the details.
 Recall that in our setting a recursive predictor is a function with range $\{0,1\}$, defined for all finite binary sequences. Axioms 1, 2, and 4 constructively define new predictors using combinations of recursive predictors. Moreover, each new predictor is defined for all inputs. Thus any new predictors constructed via the Axioms 1, 2 or 4 will also be recursive. For the partially undefined predictors obtained from Axiom 3 it suffices to specify the values of any recursive predictor in the undefined positions in order to obtain a recursive predictor satisfying Axiom 3. Thus the set of recursive predictors is closed under the axioms and therefore unpredictability with respect to this class
of predictors satisfies the universal relationship discussed in Section 5.
\hfill
$\blacksquare$

{\bf D: Alexei, would you check this proof pls?}

\section{Appendix 3: Unpredictability for different predictor classes and different predictor hierarchies}\label{appendix 3}

{\bf D: Alexei, please check the proof of Theorem 8.1. The second theorem is ok.}

Here we prove two properties of the unpredictability \eqref{IaF}.

\label{s:classes}
\begin{theorem}
There exists a non-trivial sequence $a$ with different $I(a;\F)$ for different classes $\F$ of predictors.
\end{theorem}

Proof:
Suppose we have two predictor classes, $\F = \bigcup_m$ and $\F'= \bigcup_m \F_m'. $ For predictor class $\F$, use Lemma \ref{l:blocks} to form a block $a^1$ of length $N$  which has $I(a^1:1;N)> \frac{1}{2} - \epsilon$. Form a sequence consisting of ten repeating $a^1$ blocks. Then for complexity class $\F'$, use Lemma \ref{l:blocks} to form a block $a^2$ with $I'(a^2) > \frac{1}{2} - \frac{\epsilon}{2}$. Form a sequence of $10^2$ repeated $a^2$ blocks. Continue this process to form the sequence
$$
\underbrace{a^1\ldots a^1}_{10^1\mbox{ times}}\overbrace{b^2 \ldots b^2}^{10^2\mbox{ times}}
\underbrace{a^3\ldots a^3}_{10^3\mbox{ times}}\overbrace{b^4 \ldots b^4}^{10^4\mbox{ times}}\ldots
$$
Consider the block $a^m$ with $I(a^m)> \frac{1}{2} - \frac{\epsilon}{m}$.
At the end of this block, the predicting finite state machine may be in any state. However, the class $\F_m$ consists
of all finite state machines with less than $k$ states, for some $k \in \N$. Thus finite state machines which differ
only by their starting states are all in $\F_m$. Hence $I(a^ma^m\ldots)>\frac{1}{2} - \frac{\epsilon}{m}$. Thus for
the sequence constructed above, $I(a) = \frac{1}{2}$.

We now construct Turing machine representation of a recursive predictor, and demonstrate that on the above sequence it
achieves $I(a)  = 0$. Form a tape which records the shortest repeating sequence. Use this as output. As soon as we make
a wrong prediction, find the next repeating sequence. With this machine, (guarantee a finite number of states) we will
predict perfectly somewhere in the second block, from then on, we will continue to predict perfectly until we move to
$a^{m+1}$. As soon as we accumulate errors begin to search again for the new sequence.
\hfill
$\blacksquare$

\begin{lemma}
Unpredictability is independent of the choice of hierarchy used.
\end{lemma}
Proof: Suppose we have two hierarchies of finite sets such that $\F
= \bigcup_m\F_m$ and $\F = \bigcup_m \F_m'$. Then $I(a,m)$ is
bounded below and monotonically decreasing in $m$ for both
hierarchies. We adopt the notation \eqref{IaF}, \eqref{Imn},
\eqref{Im} for the definition of the unpredictability based on the
hierarchy $\F_m$ and a similar notation $I'(a)$, $I'(a;m)$,
$I'(a;m,n)$ for the definition of unpredictability based on the
hierarchy $\F'_m$. Now, $ \F_i \subset\F= \bigcup_m \F_m' $ for each
$i$. Hence, as sets in a hierarchy are finite and increasing, there
exists a $j$ such that $ \F_i \subseteq \F_j'. $ Thus we know that
for any $i$ there exists a $j=j(i)$ such that $ I(a;i,n)\geq
I'(a;j,n)$ for all $n$. Therefore $I(a;i)\geq I'(a;j)$ and
consequently
\begin{eqnarray}
I(a)=\inf_i I(a;i)\geq \inf_j I'(a;j)=I'(a).
\end{eqnarray}
Analogously, $I'(a)\geq I(a)$. Thus $I(a) = I'(a)$. \hfill
$\blacksquare$

\subsection*{Acknowledgments}
The authors thank V. Vovk for a useful discussion of the results.
This publication has emanated from research conducted with the
financial support of Science Foundation Ireland (grant
06/RFP/MAT048) and Russian Foundation for Basic Research (grants
06-01-72552 and 06-01-00256).

{\bf D: Do we need more references?}

\bibliographystyle{beta}

\end{document}